\documentclass{amsart}

\theoremstyle{plain}

\newtheorem{thm}{Theorem}[section]

\newtheorem{prop}[thm]{Proposition}

\newtheorem{lem}[thm]{Lemma}

\theoremstyle{definition}

\newtheorem{defn}{Definition}[section]

\theoremstyle{remark}

\newtheorem{rem}{Remark}[section]

\newcommand{\Rb}{\textrm{\bf R}\xspace}

\newcommand{\pscm}[1]{\ensuremath{\mathcal R^+(#1)}}
\newcommand{\pscmt}[1]{\ensuremath{\mathcal R^+_0(#1)}}

\newcommand{\T}[1]{\ensuremath{T_{#1}}}
\newcommand{\UI}{\ensuremath{[0,1]}}
\newcommand{\GL}{\ensuremath{\mathbf{GL}}}

\newcommand{\gra}{\ensuremath{\alpha}\xspace}

\newcommand{\grb}{\ensuremath{\beta}\xspace}
\newcommand{\grG}{\ensuremath{\Gamma}\xspace}
\newcommand{\grg}{\ensuremath{\gamma}\xspace}

\newcommand{\gre}{\ensuremath{\epsilon}\xspace}

\newcommand{\grf}{\ensuremath{\theta}\xspace}

\newcommand{\grt}{\ensuremath{\tau}\xspace}

\newcommand{\grk}{\ensuremath{\kappa}\xspace}

\newcommand{\grl}{\ensuremath{\lambda}\xspace}

\newcommand{\grp}{\ensuremath{\pi}\xspace}

\newcommand{\grd}{\ensuremath{\delta}\xspace}

\newcommand{\grx}{\ensuremath{\xi}\xspace}
\newcommand{\grPh}{\ensuremath{\Phi}\xspace}
\newcommand{\grph}{\ensuremath{\phi}\xspace}

\newcommand{\grs}{\ensuremath{\sigma}\xspace}
\newcommand{\grPs}{\ensuremath{\Psi}\xspace}
\newcommand{\grps}{\ensuremath{\psi}\xspace}

\newcommand{\grn}{\ensuremath{\nu}\xspace}

\newcommand{\grm}{\ensuremath{\mu}\xspace}

\newcommand{\grr}{\ensuremath{\rho}\xspace}
\newcommand{\grUpsilon}{\ensuremath{\Upsilon}\xspace}

\usepackage{latexsym,amsmath}
\usepackage[dvips]{graphics}
\usepackage{amssymb,xspace,amscd}
\usepackage[all]{xy}

\begin{document}

% Top matter stuff

\title[On the homotopy type of $\mathcal{R}^+(M)$]{On the homotopy type of the space $\mathcal{R}^+(M)$}
\author{Vladislav Chernysh}
\address{ Department of Mathematics\\ University of Notre
Dame\\Notre Dame, IN 46556}%
\email{vchernys@nd.edu}%
\subjclass[2000]{58D17, 57R65}%

\begin{abstract} In this paper it is shown that the space of metrics of
positive scalar curvature on a manifold is, when nonempty,
homotopy equivalent to a space of metrics of positive scalar
curvature that restrict to a fixed metric near a given submanifold
of codimension greater or equal than $3$. Our main tool is a
parameterized version of the Gromov-Lawson construction, which was
used to show that the existence of a metric of positive scalar
curvature on a manifold was invariant under surgeries in
codimension greater or equal than $3$.
\end{abstract}

\maketitle

% End top matter stuff

\section{Introduction}\label{intro}

The topology of the space of metrics of negative scalar curvature
\ensuremath{\mathcal R^-(M^n)} was described by Lohkamp
\cite{Loh92}, and it turns out that \ensuremath{\mathcal R^-(M^n)}
is nonempty and contractible, when $n\ge 3$ and $M^n$ is closed.

On the other hand, the questions of existence and classification
for metrics of positive scalar curvature are much more
complicated. There are manifolds that do not admit a metric of
positive scalar curvature. Any spin manifold $M^n$ for which $n\ge
5$ and $\hat{A}(M^n)\ne 0$ is such a manifold. On the
classification side, R.~Carr~\cite{Car89} shows that the space of
positive scalar curvature metrics on a sphere $S^{4m-1}$, $m\ge 2$
has infinitely many connected components. Also, in general, it is
not known when a given closed manifold admits a metric of positive
scalar curvature.

In this paper, we concern ourselves with the homotopy type of the
space $\pscm{M}$. Since Surgery Theory is a major tool in studying
manifolds, it is important to understand how the homotopy type of
$\pscm{M}$ changes under surgeries.

Let $M^n$ be a closed manifold and $\pscm{M^n}$ be the space of
metrics of positive scalar curvature on $M^n$ (which is assumed to
be nonempty throughout this paper). The topology on this space is
defined by the collection of $C^k$-norms $\|.\|_k$ on the space of
all Riemannian metrics $\mathcal R (M^n)$ with respect to some
reference metric $h$: $\|g\|_k=\max_{i\le k} \sup_{M^n}|\nabla^i
g|$. This topology does not depend on the choice of the metric
$h$. Let $N^{n-k}\subset M^n$, $k\ge 3$, be a closed submanifold
with a trivial normal bundle. We fix a tubular neighborhood
$\grt\colon N\times D^k\rightarrow M$, an arbitrary metric $g_N$
on $N$, and a torpedo metric $g_0$ on $D^k$ (for the definition of
a torpedo metric see below, page~\pageref{torpedo}), such that the
metric $g_N+g_0$ has positive scalar curvature on $N\times D^k$.

Given these data, we define
\begin{equation}\label{Ro}
\pscmt{M^n}:=\{g\in\pscm{M^n}|\grt^*(g)=g_N+g_0\}.
\end{equation}
Our main result is the following theorem.
\begin{thm}\label{main} Suppose that $\pscm{M^n}$ is not empty. Then
the inclusion map \[ i\colon \pscmt{M}\rightarrow\pscm{M^n}\] is a
homotopy equivalence.
\end{thm}
As an easy corollary we have the Surgery Theorem.
\begin{thm}\label{surgery}
Let $S^k\rightarrow M_1^n$ be an embedding with the trivial normal
bundle and $\grt_1\colon S^k\times D^{n-k}\rightarrow M_1^n$ be a
tubular neighborhood for this embedding. Let $M_2$ be a manifold
obtained by doing surgery on this tubular neighborhood. If $n-k\ge
3$ and $k\ge 2$ then $\pscm{M_1}$ and $\pscm{M_2}$ are homotopy
equivalent.
\end{thm}

The idea of the proof is as follows. From Palais~\cite{Pal66}, it
follows that $\pscm{M}$ and $\pscmt{M}$ are dominated by
$CW$-complexes. Therefore, it suffices to show that the homotopy
sets $\grp_k\left(\pscm{M},\pscmt{M}\right)=0$, all $k\ge 0$, and
the inclusion map $i$ is a bijection on path components of both
spaces.

The Gromov-Lawson construction gives us a deformation $\GL$ of a
compact family $g_s\in\pscm{M}$ into $\pscmt{M}$. The only problem
is that $\pscmt{M}$ is not invariant under this deformation. The
saving grace of the Gromov-Lawson construction is its local
$O(k)$-symmetry on normal fibres of $N$.

Thus, we overcome this problem by introducing the space
$W(N,\grt)$ of metrics with a warped fibre. The space $W(N,\grt)$
is invariant under $\GL$ and we show that the subspace $\pscmt{M}$
is a weak deformation retract of $W(N,\grt)$.

The Greek letter $\grk$ is used to denote the scalar curvature
throughout the paper.

I would like to thank Stephan Stolz, Bill Dwyer, and Larry Taylor
for numerous discussions and their valuable advice with regard to
this work.

\section{The Gromov-Lawson construction}

The Gromov-Lawson construction (see \cite{GL80}, \cite{RS01})
allows one to conclude that if a closed manifold $M_2$ is obtained
from a closed manifold $M_1$ by doing surgery in codimension $\ge
3$ and $M_1$ carries a metric of positive scalar curvature, then
$M_2$ carries a metric of positive scalar curvature. It is our
main tool in deforming a compact family $g_s\in\pscm{M}$ into
$\pscmt{M}$, see Theorems~\ref{isolemma},~\ref{alphatwo},
and~\ref{alphathree}.

\label{ho}Throughout this paper, we fix a tubular neighborhood
$\grt\colon N\times D^k_{T_0}\rightarrow M$, an arbitrary metric
$g_N$ on $N$, and a torpedo metric $g_0$ on
$D^k_{T_0}\subset\Rb^k$ of radius $\gre_0$ and size $T_0=\gre_0T$
(see the definition below), such that $g_N+g_0$ has positive
scalar curvature on $N\times D^k$. Since we assume that $\pscm{M}$
is not empty, the original Gromov-Lawson construction gives us a
metric $h_0$ on $M$ (of positive scalar curvature), which is
isometric to $g_N+g_0$ near $N$. We also fix a vector bundle
isometry $\grph\colon (N\times\Rb^k, g_{\mathrm{eucl}})\rightarrow
(\perp N,h_0)$ (recall that one of our assumptions is that the
normal bundle of $N$ is trivial). We take the restriction
$\grph\colon N\times D_{T_0}\rightarrow \perp N$ and define
$\grt_0:=\exp ^{\perp}_{h_0}\circ\,\grph\colon N\times
D^k_{T_0}\rightarrow M$, where $\exp^{\perp}_{h_0}$ is the normal
exponential map for the metric $h_0$. Then
$\grt_0^*(h_0)=g_N+g_0$. By the uniqueness theorem for tubular
neighborhoods we may assume that $\grt=\grt_0$.

\begin{figure}\label{figintro} \centering
\includegraphics{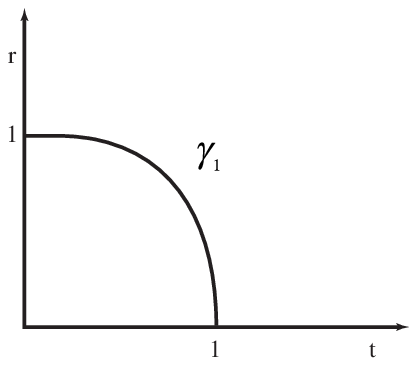}
\caption{}\label{Figure8}
\end{figure}

In general, by a torpedo metric in the disc $D_{T_0}^k$, we
understand an $O(k)$-symmetric positive scalar curvature metric
which is a Riemannian product with a standard euclidean $k-1$
sphere near the boundary and is a standard euclidean $k$ sphere
metric near the center of the disc. More specifically, we fix a
curve $\grg_1$ in the plane $(t,r)$ as in Figure~\ref{Figure8},
which follows a horizontal line near the point $(0,1)$, and then
joins smoothly an arc of a circle of radius
$1$.The curve $\grg_{\gre}$ is the
image of $\grg_1$ under the homothety of the plane with the
coefficient $\gre$. If $T$ is the length of $\grg_1$, then $\gre
T$ is the length of $\grg_{\gre}$. Let
\begin{equation*}\label{torpedo}
T_{\grg_{\gre}}:=\{(x,t)\in \Rb^{k}\times\Rb\ |\
(t,|x|)\in\grg_{\gre}\}.
\end{equation*}
Here the metrics on $\Rb^{k}$ and $\Rb$ are the standard euclidean
metric. We take $T_{\grg_{\gre}}$ together with the induced metric
and call such a metric a torpedo metric of radius $\gre$ and size
$\gre T$. \suppressfloats

\subsection{The Gromov-Lawson construction}

We begin by briefly describing the main construction by following
\cite{RS01}. For details, see \cite{GL80}, \cite{RS01}.

We take a Riemannian product $M\times \Rb$ and consider the
restriction of the product metric to the hypersurface

\begin{equation}\label{neck}
\T{\grg}(g)=\{(y,t)\in M\times\Rb \,|\, (t,\|y\|)\in \grg\},
\end{equation}

\noindent where \grg\ is a  $C^{\infty}$-curve of finite length
which lies in the first (NE) quadrant of $\Rb^2$ plane with
coordinates $(t,r)$, see Figure~\ref{Figure1}, and $\|y\|$ is the
distance from $y$ to $N$ with respect to the metric $g$. We call
this hypersurface a neck.

\begin{defn} Let $g_s$, $s\in S$ be a family of pscm metrics
continuously parameterized by a compact space $S$ and
$\grg\colon\Rb\rightarrow\Rb^2=\{(t,r)\}$ be a $C^{\infty}$
isometric embedding. We call $\grg$ \emph{an admissible curve for
a family} $g_s$ if the following holds:

\begin{enumerate}
\item $\grg(0)=(0,r_0)$ with $r_0>0$ is such that $\grg$ follows
the $r$-axis linearly from $\grg(-\infty)=(0,\infty)$ to
$(0,r_0)$;

\item $\grg$ intercepts the $t$-axis at a right angle and follows
an arc of a circle (of possibly infinite radius) at this point;

\item the curve $\grg$ crosses the line $r=r_0$ only once and is
symmetric about the $t$-axis, i.e. $\grg(L-s)=R_t\circ\grg(L+s)$,
where $L$ is a unique number such that $\grg(L)\in
t\mathrm{-axis}$ and $R_t$ is the reflection about $t$-axis;

\item the injectivity radius of the normal exponential map for all
$g_s$ is strictly greater than $r_0$.
\end{enumerate}

\end{defn}

The space of all admissible curves is denoted by $\grG$ (note that
the vertical segment is an admissible curve) and has a natural
topology which is induced from the $C^{\infty}$-topology on
$C^{\infty}(\Rb,\Rb^2)$.

\begin{rem}
(i) Any curve $\grg\in\grG$ is uniquely determined by its part on
$[0,L]$; (ii) For an admissible curve, the neck $T_{\grg}$ can be
defined over the tubular neighborhood of $N$ of radius $r_0$ by
taking the part of the curve on $[0,L]$ and applying the
formula~\ref{neck} above.
\end{rem}

\begin{prop}\label{emb}
Let $\{g_s\}$ be a family of metrics parameterized by a compact
space $S$ and $\grG$ be the space of admissible curves for this
family. Then there exists a continuous map\[
S\times\grG\rightarrow \mathrm{Emb}^{\infty}(M,M\times\Rb)\] with
the following properties: (i) for any pair $(s,\grg)$ the
resulting embedding $f$ is constant outside of the normal tubular
neighborhood $Tb_{r_0}(N)$ of radius $r_0$ in the sense that
$f(x)=(x,0)$ for all $x\in M-Tb_{r_0}(N)$; (ii) $f$
diffeomorphically maps $Tb_{r_0}(N)$ onto the neck $T_{\grg}$.
\end{prop}

\begin{proof}

We fix a family of smooth increasing functions $\grph$:
$\grph(r)=r+r_0-L$ if $r\ge L\ge r_0$, $\grph(r)=r$ if $r\le\grd$,
some $0\le\grd<r_0$. Let $N_{\grg}$ be a submanifold of $T_{\grg}$
which is diffeomorphic to $N$ under the projection map $p\colon
M\times\Rb\longrightarrow M$ and $g_{\grg}$ is a metric on the
neck induced from the product metric on $M\times\Rb$.

Recall that $L$ is the length of the curve $\grg$ between the
point where it intercepts the $t$-axis and the point $(0,r_0)$.
The function $\grph$ defines a diffeomorphism of the normal bundle
$\grph\colon\grn(N)\rightarrow\grn(N)$,
$(x,v)\mapsto(x,\frac{\grph(|v|)}{|v|}v)$. Then the composition of
maps
\begin{eqnarray*}
\lefteqn{Tb_{L}(N_{\grg})\xrightarrow{(\exp^{\perp}_{g_{\grg}})^{-1}}\grn(N_{\grg})\xrightarrow{d p}%
\grn(N)\xrightarrow{\grph}}\\
& & \grn(N)\xrightarrow{\exp^{\perp}_g}Tb_{r_0}(N)
\end{eqnarray*}
is a diffeomorphism of tubular neighborhoods which is the identity
on all points whose distance from $N$ is greater or equal than
$r_0$ and less than $r_0+\gre''$, for some small enough
$\gre''>0$; all distances are taken with respect to the metric
$g$.

It is clear that we may add the points $x\in M$ outside this
neighborhood as $(x,0)$. Thus we obtain an embedding of $M$ into
$M\times\Rb$.
\end{proof}

The induced metric on the neck can be considered as a metric on
the manifold $M$ via the pullback by the embedding that we obtain
from the lemma above.

We call an admissible curve $\grg$ horizontal if (i) the
coordinates $r(s)$, $t(s)$ are monotone functions of $s$ when
$s\in(-\infty,0]$; (ii) it follows a segment parallel to $t$-axis
somewhere between $r=r_0$ and $r=0$.

\begin{thm}\label{curve}
Let $N$ be a closed submanifold of $M$ of codimension $k\ge 3$ and
$g_s$ be a continuous family of metrics of positive scalar
curvature, parameterized by a compact space $S$. Then there exist
a number of points $(t_i,r_i)\in\Rb^2$, $i=0..5$, see
Figure~\ref{Figure1}, and a horizontal (along $[t_4,t_5]\times
r_4$) admissible curve $\grg$ such that the scalar curvature of
the neck $T_{\grg}(g_s)$ is positive for all $s\in S$.
\end{thm}

\begin{proof}

\begin{figure}

\centering
\includegraphics{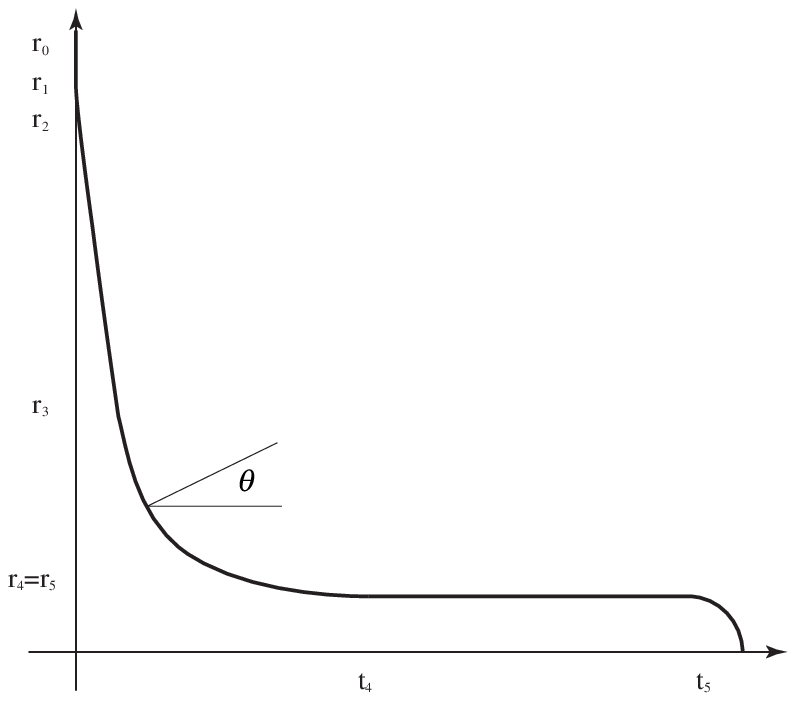}
\caption{}\label{Figure1}

\end{figure}

This is essentially proved in the original work~\cite{GL80} and
improved in~\cite{RS01}. We only add that the argument goes
through for a \emph{compact} family of metrics.

The first bend is small and occurs between $r_1$ and $r_2$ and is
followed by a straight line. The second bend occurs between $r_3$
and $r_4$ and ends in a horizontal segment. If $\mathbf{k}$ is the
signed curvature of the curve $\grg$, then during the second bend
$\mathbf{k}\le\frac{\sin\grf}{2r}$. Here, $\grf$ is the angle
between normal to $\grg$ and $t$-axis. In addition, we can arrange
that
\begin{eqnarray} \mathbf{k}\le\frac{9\sin\grf}{16r}\label{curvestimate}
\end{eqnarray}
will still assure the positivity of the scalar curvature during
the second bend.

From the equations controlling the scalar curvature on the neck we
see that the scalar curvature will always be positive on the part
bending towards the $t$-axis ($\mathbf{k}\le 0$) provided $r_4$ is
small enough. We assume that the end part of $\grg$ is exactly the
torpedo curve $\grg_{\gre}$ (for $\gre=r_4$).
\end{proof}

\subsection{Initial bend}

Let $\grg$ be a curve in the plane $(t,r)$, as in the statement of
Theorem~\ref{curve}. In this section we prove that one can deform
$\grg$ through admissible curves into the $r$-axis, keeping the
scalar curvature on the neck positive.

Any admissible curve is uniquely determined by its curvature
function $\mathbf{k}(s)$ on the interval $[0,L]$. So, we deform
such a curve by deforming its curvature function.

Let $\grd_{s_0}$ be a \grd-cutoff function at a point
$s_0\in[s_3,s_4]$, i.e. a $C^{\infty}$-function equal to $1$ for
all $s\le s_0$, and equal to $0$ for all $s \ge s_0+\grd$ (the
interval $[s_3, s_4]$ corresponds to the part of the curve between
$r_3$ and $r_4$). Let
$\tilde{\grg}:=\grg(\tilde{\mathbf{k}}_{s_0})$ be the curve
corresponding to the curvature function
$\tilde{\mathbf{k}}_{s_0}(s):= \grd_{s_0}(s)\mathbf{k}(s)$.

\begin{prop}\label{kprop} There exists a $\grd>0$ such that for any
$s_0\in [s_3,s_4]$ the curvature $\tilde{\mathbf{k}}$ (which is
equal to $\tilde{\mathbf{k}}_{s_0}$) of the curve $\tilde{\grg}$
satisfies
\[
\tilde{\mathbf{k}}(s)\le\frac{9}{16}\frac{\sin\tilde{\grf}(s)}{\tilde{r}(s)}
\]
on this interval. Here, $\tilde{\grf}$ is the angle between the
normal to the curve $\tilde{\grg}$ and the $t$-axis.
\end{prop}

\begin{proof}
Since $\grg$ is the curve from Theorem~\ref{curve}, its curvature
satisfies $\mathbf{k}\le \frac{\sin \grf}{2r}$ on the interval
$[s_3,s_4]$. Consider the function
\begin{eqnarray*}
F\colon[s_3,s_4]\times[0,1]&\rightarrow& \Rb\\
(s,t)&\mapsto&
\mathbf{k}(s+t)-\frac{9}{16}\frac{\sin\grf(s)}{r(s)}.
\end{eqnarray*}
We have that $F([s_3,s_4]\times 0)\le A_0 <0$, so there exists a
$\grd>0$ such that $F([s_3,s_4]\times[0,\grd])\le A_1<0$. We take
a cutoff function corresponding to \grd and observe that for all
$0\le s \le \grd$ and all $s_0\in[s_3,s_4]$
\begin{eqnarray*}
\tilde{\mathbf{k}}(s_0+s)&=&\grd_{s_0}(s_0+s)\mathbf{k}(s_0+s)\\
&\le& \mathbf{k}(s_0+s)<\frac{9}{16}\frac{\sin \grf(s_0)}{r(s_0)}\\
&\le&
\frac{9}{16}\frac{\sin\tilde{\grf}(s_0+s)}{\tilde{r}(s_0+s)}.
\end{eqnarray*}
\end{proof}

\par\noindent Now we can prove the following theorem.

\begin{thm}\label{isolemma} Suppose $\grg\in\grG$ is as in Theorem~\ref{curve}, where $\grG$ is the
set of admissible curves for the family $g_s$ .Then there exists a
continuous map
$\gra_1:[0,1]\to\grG$ such that%
\begin{enumerate}
\item $\gra_1(1)=\{r\textrm{-axis}\}$, $\gra_1(0)=\grg $;%
\item the scalar curvature on each $\T{\gra_1(t)}$ is positive for
all $g_s$.
\end{enumerate}
\end{thm}

\begin{proof}

\begin{figure}

\centering
\includegraphics{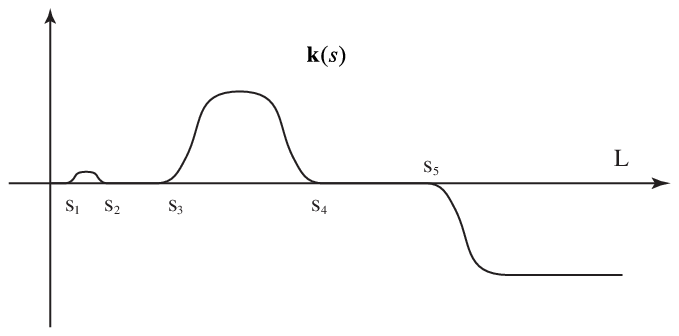}
\caption{}\label{Figure2}

\end{figure}

From the proof of Theorem~\ref{curve} we may assume that the
signed curvature of $\grg$ is given by a $C^{\infty}$- function
$\mathbf{k}\colon [0, L] \rightarrow \Rb$ with the graph as in
Figure~\ref{Figure2}, and $\int_0^L \mathbf{k}(s)ds=0 $. In this
picture $s_4=\sup\{s\ge 0 | f(s)>0\}$ and $s_5=\inf\{s\ge 0 |
f(s)<0\} $ are positive numbers that remain fixed throughout the
argument.

We deform the curve by deforming its curvature function to $0$.
The delicate part is to keep the first and the second bends where
they initially took place, and to keep the last bend confined
within the region $0\le r\le r_4$.

Now, we choose a small number $0< \gre <\min(s_1,s_5-s_4)$ and a
linear path $p(t)= s_4 - s_4*t$. We take \gre to be equal to \grd
of Proposition~\ref{kprop} and take a cutoff function
$\gre_{p(t)}$ at $p(t)$ with the parameter \gre. The deformation
function $\mathbf{k}_t$ is $\tilde{\mathbf{k}}_{p(t)}$ on
$[0,p(t)+\gre]$, and a carefully rescaled function on
$[p(t)+\gre,L_t]$ where it is non-positive. The curve
$\grg(\tilde{\mathbf{k}}_{p(t)})$ will satisfy the
inequality~\ref{curvestimate} and therefore will give the positive
scalar curvature on the neck.

We set $\gra_1(t):=\grg(\mathbf{k}_t)$.\end{proof}

From Palais~\cite{Pal68} and Proposition~\ref{emb}, we conclude that the pullback map %
\[\mathrm{Emb}^{\infty}(M,M\times\Rb)\times\mathcal{R}(M\times\Rb)\rightarrow\mathcal{R}(M)\]
is continuous. Therefore, we may regard $\gra_1$ as a continuous
deformation of a family of metrics $g_s$ on $M$.

\subsection{Middle stage deformation}

During the second step of our construction, we further deform a
given metric to a metric which on some tubular neighborhood of $N$
is a Riemannian product of $g_N$ and a torpedo metric of some
small radius $r_4$.

Given a smooth family of metrics $\gra\colon I\rightarrow
\pscm{X}$ on a closed smooth manifold $X$, we can put a metric on
$X\times\Rb$. However, in general, the scalar curvature of this
obvious metric $g(x,t)=\gra(t)(x)+dt^2$ will not be positive.

We fix a function $F\colon\Rb\rightarrow [0,1]$ such that $0\le
F'<2$, $F(t)=0,\, t\in(-\infty,\gre]$, $F(t)=1,\,
t\in[1-\gre,\infty)$, for some $0<\gre<1/4$, and for a positive
number $\grt$ we define a function $F_{\grt}\colon\Rb\rightarrow
\Rb$, by $F_{\grt}(t)=F(\frac {t}{\grt})$.

Let
\begin{equation}\label{product}
g^{\gra}_{\grt}(x,t):=\gra(F_{\grt}(t))(x)+dt^2.
\end{equation}

By Lemma on page $184$ in~\cite{Gaj87}, there exists a number
$\grt>0$ such that this metric has positive scalar curvature and
is a Riemannian product near $X\times 0$ and $X\times t_0$.

\begin{prop}\label{square}
Let $\grn=(E, B, p)$ be a finite dimensional vector bundle, and
$g_1$, $g_2$ be two bundle metrics on $\grn$. Then there exists a
canonical vector bundle isometry:
\begin{equation*}
\xymatrix{ (E,g_1) \ar[rr]^{\sqrt{(g_1,g_2)}} \ar[dr]_p & & (E,g_2) \ar[dl]^p\\
& B} .\end{equation*}%
\end{prop}

\begin{proof}
For each fiber $V$ there exists a unique positive symmetric
operator $A$ on $V$, such that $g_1(u,v)=g_2(Au,v)$. From the
spectral theorem, there exists a unique positive symmetric
operator $B$ with the property $B^2=A$. We have:
$g_1(u,v)=g_2(Au,v)=g_2(BBu,v)=g_2(Bu,Bv)$. The operator
$\sqrt{(g_1,g_2)}:=B\colon V\rightarrow V$ provides the required
isometry.
\end{proof}

Let $\grg$ be the curve from Theorem~\ref{curve}. By $\tilde{g}_s$
we denote the family of metrics obtained by pulling back the
induced metrics on the necks $T_{\grg}(g_s)\subset M\times\Rb$ via
our canonical diffeomorphisms.

Our goal is to deform the family $\tilde{g}_s$ further, so that
the resulting metrics near $N$ are isometric to $g_N+g_0$, where
$g_0$ is a torpedo metric of radius $1$.

The key observation here is the following. If we start with an
arbitrary metric on $M$ (not necessarily of positive scalar
curvature) then the metric induced on the $\gre$-sphere bundle of
$N$ will have positive scalar curvature if $\gre$ is sufficiently
small. Here, the assumption that codimension is $\ge 3$ is
crucial. In the same fashion, if $T_{\grg_{\gre}}$ is a ``cap"
(i.e. the neck for some torpedo metric curve $\grg_{\gre}$) of
some small radius $\gre$ over a neighborhood of $N$, then the
scalar curvature of $T_{\grg_{\gre}}$ is positive.

\begin{thm}\label{alphatwo}
There exists a continuous map
\begin{equation*}
\gra_2\colon S\times I\rightarrow\pscm{M}
\end{equation*}
such that: (i) $\gra_2(s,0)=\tilde{g}_s$; (ii)
$\grt^*_s(\gra_2(s,1))=g_N+g_0$, where $g_N$ is a fixed metric on
$N$, $g_0$ is the fixed torpedo metric of radius $\gre_0$, and
$\grt_s:=\exp^{\perp}_{\gra_2(s,1)}\circ\sqrt{(h_0,\gra_2(s,1))}\circ\grph\colon
N\times D^k_{T_0}\rightarrow M$ is a tubular map for the metric
$\gra_2(s,1)$.
\end{thm}

\begin{proof}

Given a compact family of metrics $\tilde{g}_s$ on $M$, we define
a family of metrics on the tubular neighborhood $\grt(N)$:
\begin{eqnarray*}
\bar{\gra}\colon S\times I&\rightarrow&\mathcal{M}(\grt(N)),\\
(s,t)&\mapsto& (1-t)\tilde{g}_s+t(g_N+g_{\mathrm{eucl}}).
\end{eqnarray*}
Here, $\tilde{g}_s$ denotes the restriction of $\tilde{g}_s$ to
$\grt(N)$, the metric $g_{\mathrm{eucl}}$ is the standard flat
metric on $D^k$ and the product is taken with respect to our fixed
tubular map $\grt\colon N\times D^k\rightarrow M$. We take a curve
$\grg_{\gre}$ corresponding to a torpedo metric of radius $\gre$
and a cap
\[T_{\grg_{\gre}}(\bar{\gra}(s,t)):=\{(y,t)\ |\
(t,\|y\|)\in\grg_{\gre}\} \] in $\grt(N)\times\Rb$. The distance
$\|y\|$ from $y$ to $N$ is taken with respect to
$\bar{\gra}(s,t)$. We choose $\gre$ such that for all $(s,t)\in
S\times I$ the metric $\bar{\gra}(s,t)$ is a diffeomorphism from
the normal disc bundle of radius $\gre$ into $\grt(N)$, and the
scalar curvature of the induced metric on the cap
$T_{\grg_{\gre}}(\bar{\gra}(s,t))$ is positive. For such a choice
of $\gre$ all caps are canonically diffeomorphic via normal
exponential maps and the canonical isometry
$\sqrt{(g_1,g_2)}\colon(\perp N, g_1)\rightarrow(\perp N, g_1)$
for any two such metrics $g_1$ and $g_2$.

Let $\grg\in$ be the curve in the plane $(t,r)$, constructed in
Theorem~\ref{curve}. Recall that this curve comes with a choice of
a number of parameters $r_i$, $t_i$. In particular, when $r=r_4$,
it follows the horizontal line between $t_4$ and $t_5$. Without
loss of generality we may assume that $\gre=r_4$.

Let $s\in S$ and $T_{\grg}\subset M\times\Rb$ be the neck
corresponding to $\grg$ and the metric $g_s$. Let $N'$ be the part
of the neck for $t\ge t_5$, $N''$ the part for $t_4\le t\le t_5$,
and $N'''$ the part for $t\le t_4$. Recall that the part of the
curve that defines $N'$ is exactly the curve $\grg_{\gre}$. The
boundary of $N'$ is the spherical bundle $S(N)$ of $N$ of radius
$\gre$, which we consider as a submanifold of $M$.

We pull the induced metric on the cap for $g_t:=\bar{\gra}(s,t)$
to $N'$ via our canonical diffeomorphism. This gives us a smooth
family of metrics on $N'$, $\grb_1\colon I\rightarrow \pscm{N'}$.
The metric $\grb_1(1)$ is isometric to $g_N+g_{\gre}$, where
$g_{\gre}$ is the torpedo metric of radius $\gre$. There is a
canonical path $\grb_2$ from $g_{\gre}$ to our fixed torpedo
metric $g_{0}$. It is obtained by taking a linear path from $\gre$
to $\gre_0$. Applying this path to fibers of the normal bundle we
obtain a path $\grb\colon I\rightarrow\pscm{N'}$ as follows
\begin{equation*}
\grb(\grt):=\begin{cases} \grb_1(2\grt),& \grt\in[0,1/2]\\
\grb_2(2(\grt-1/2)),& \grt\in[1/2,1]\end{cases}.
\end{equation*}
We have, $\grb(1)=g_N+g_0$. The path constructed in such a way is
not necessarily smooth in $\grt$ at the point $\grt=1/2$. However,
we can always apply the standard smoothing procedure for piecewise
smooth paths. So, we assume $\grb$ is smooth in $\grt$ everywhere.

By the formula~\ref{product}, we obtain a metric
$g(x,t):=\bar{\grb}(F_{t_0}(t-t_4))(x)+dt^2$ of positive scalar
curvature on $S(N)\times[t_4,t_5]$, where
$t_0=\frac{1}{(t_5-t_4)}$. In this formula, $\bar{\grb}$ is the
restriction of $\grb$ to $S(N)$. For a fixed $t_4$ we can always
choose $t_5$ large enough, so that $t_0$ is small enough to assure
positivity of the scalar curvature of $g(x,t)$. Near the boundary
$S(N)\times t_4$ this metric is a Riemannian product of the
restriction of $g_s$ to $S(N)$ and the standard metric on $\Rb$,
near the boundary $S(N)\times t_5$ it is a Riemannian product of
the metric $\grb(1)|_{\partial N'}$ and the standard metric on
$\Rb$.

We now define a map
\begin{eqnarray}
\gra_2&\colon& S\times I \rightarrow\pscm{M}\\
\gra_2(s,\grt)&=&\begin{cases}g_s|_{T_{\grg}},& \textrm{ on }
N'''\\
\bar{\grb}(\grt F_{t_0}(t-t_4))+dt^2,& \textrm{ on }N''\\
\grb(\grt),& \textrm{ on }N'
\end{cases}\nonumber
\end{eqnarray}

This map is continuous and it gives us the required
deformation.\end{proof}

\subsection{Final deformation}

After completing the second deformation, all metrics in the new
family $\gra_2(s,1)$, $s\in S$ have the required form $g_N+g_0$
near the submanifold $N$ with respect to their individual tubular
maps. However, these tubular maps are not necessarily the same as
our fixed tubular map $\grt$. Therefore, we need an explicit
version of the unique tubular neighborhood theorem for a family of
metrics.

Let $\{ g_s\}$ a family of metrics on $M$ parameterized by a
smooth compact manifold $S$. For our fixed metric $h_0$ we define
a new family of metrics $\{ g_{s,t}:= (1-t)h_0 + tg_s\}_{s\in S,\,
t\in\UI}$ on $M$.

For any two numbers $\gre^*,\; \gre^{**}>0$ we  fix a radial
diffeomorphism $\grps\colon\Rb^k\rightarrow \Rb^k$ with the
properties:
\begin{itemize}
\item[(i)] $\grps(\Rb^k)\subset
D_{\gre^*+\gre^{**}}\Rb^k$%
\item[(ii)] $\grps|_{D_{\gre^*}}=\mathrm{id}$.
\end{itemize}Radial means $\grps(x)=\grl(|x|) x$ with
$\grl > 0$.

Since $S$ and $N$ are compact, we may find two positive numbers
$\gre^*$ and $\gre^{**}$ such that for each pair $(s,t)$ the map
$\grt_{s,t}:=
\exp^{s,t}\circ\,(\grps)\circ\sqrt{(h_0,g_{s,t})}\circ\grph%
\colon N\times\Rb^k\rightarrow M^n$ is an embedding. Here, the map
$\sqrt{(h_0,g_{s,t})}$ is the canonical isometry between $(\perp
N, h_0)$ and $(\perp N, g_{s,t})$. Taking the restriction
$h_{s,t}:=\grt_{s,t}|_{N\times D_{\gre^*}}$ for a fixed $s\in S$,
we obtain an isotopy of tubular neighborhoods of $N$, where each
neighborhood has radius $\gre^*$ with respect to the metric
$g_{s,t}$. From the parameterized version of the Thom's embedding
theorem it follows that this family of isotopies is embeddable,
i.e. there exist a compact neighborhood $V_0$ of $N$ and a
$C^0$-family $\{H_{s}\}$ of diffeotopies of $M^n$ such that
$h_{s,t}=H_{s,t}\circ\;h_{s,0}$ and all $H_{s,t}$ leave points
outside $V_0$ fixed.

We deform the family $g_s$ further, so that the resulting family
is \emph{equal} to $h_0$ on the fixed tubular neighborhood
$\grt(N)\subset M$.

\begin{thm}\label{alphathree} Let $g_s=\gra_2(s,1)$, $s\in S$ be a continuous family
of metrics parameterized by a compact manifold $S$, such that
$\grt^*_s(g_s)=g_N+g_0$. Then there exists a continuous map
\begin{equation*}
\gra_3\colon S\times I\rightarrow\pscm{M},
\end{equation*}
such that $\gra_3(s,0)=g_s$; $\grt^*(\gra_3(s,1))=g_N+g_0$ on
$N\times D^k_{T_0}$, where $\grt\colon N\times
D^k_{T_0}\rightarrow M$ is the fixed tubular map.
\end{thm}

\begin{proof}
Without loss of generality, we can assume that for each $g_s$ the
normal exponential map is a diffeomorphism on the normal disc
bundle of radius $T_1>T_0$. For a number $0<c\le 1$ and a metric
$g_s$ we define a diffeotopy of $M$ as follows. Fix a family of
increasing functions $\grph_{t}$ on $[0,T_1]$ parameterized by
$t\in[0,1]$, such that $\grph_0(\grt)=\grt$ and
$\grph_1(\grt)=c\grt$ on $[0,T_0]$, $\grph_1(\grt)=\grt$ near
$T_1$. This family defines a family $\grPh_{s,t}$ of diffeotopies
of $M$ by a diffeomorphism of a normal bundle
$(x,v)\mapsto(x,\frac{\grph_t(|v|)}{|v|}v)$ and then applying a
normal exponential map of $g_s$. Now, let $H_{s,t}$ be a family of
diffeotopies of $M$ corresponding to the family $g_{s,t}$ defined
above, and $\gre^*$ be the radius of the tubular neighborhood
which is embedded by this family of diffeotopies. We may always
choose $\gre^*$ so that $c=\frac{\gre^*}{T_0}\le 1$. We define
\begin{equation*}
\bar{\grPh}_{s,t}:=\begin{cases}\grPh_{s,3t},& t\in[0,1/3]\\
H_{s,3t-1}^{-1}\circ\grPh_{s,1},& t\in[1/3,2/3]\\
\grPh^{-1}_{0,3t-2}\circ H^{-1}_{s,1}\circ\grPh_{s,1},&t\in[2/3,1]
\end{cases},
\end{equation*}
where $\grPh_{0,t}$ is the family of diffeotopies for the fixed
metric $h_0$.

The map $\gra_3$ is defined as
\begin{equation*}
\gra_3(s,t):=\left(\bar{\grPh}^{-1}_{s,t}\right)^*(g_s).
\end{equation*}

The differential of $\bar{\grPh}_{s,1}$, when restricted to the
fibers of the normal bundle is the canonical vector bundle
isometry $\sqrt{(g_s,h_0)}$ between $(\perp N,g_s)$ and $(\perp N,
h_0)$.

We have $\grt^*(\gra_3(s,1))=\left(\grPh^{-1}_{s,1}\circ
H_{s,1}\circ\grPh_{0,1}\circ\grt\right)^*(g_s)$. An easy
computation shows that $\grPh^{-1}_{s,1}\circ
H_{s,1}\circ\grPh_{0,1}\circ\grt=\grt_s$.
\end{proof}

\section{Locally warped metrics}

We fix a nonnegative number $B$ and define a space
\[ W:=\{h\in\pscm{D^n_{T_0}}\, |\, h= g(t)^2dt^2+f(t)^2d\grx^2\}\]
of warped metrics in the disc with the additional requirement that
the scalar curvature of $h$ is greater than $B$ everywhere in the
disc. The metric $d\grx^2$ is the standard metric on the euclidean
sphere of radius $1$.

The subspace $W^{\mathrm{loc}}\subset W$ consists of metrics that
are isometric to the torpedo metric $g_0$ near the origin $0\in
D^n_{T_0}$, see Definition~\ref{loc}.

The main result of this section is Theorem~\ref{locdef}, which
states that there exists a deformation of $W$ into
$W^{\mathrm{loc}}$ such that it does not change the metric near
the boundary of the disc.

\subsection{Warped products}

Here we briefly recall some facts about warped Riemannian
products.

\begin{defn}
Let $(B,\check{g})$ and $(F,\hat{g})$ be Riemannian manifolds, and
$f$ a positive function on $B$. Then the Riemannian manifold
$(B\times F, \check{g}+f\hat{g})$ is called a warped product.
\end{defn}

In particular, any warped product is a Riemannian submersion over
the base $(B,\check{g})$ with respect to the projection $\grp$ on
the first factor. To any Riemannian submersion $\grp\colon
(M,g)\rightarrow (B,\check{g}) $ we can correspond two fundamental
invariants $A$ and $T$, which are $(2,1)$ tensors on $M$.

\begin{defn} For any two vector fields $E_1$ and $E_2$ on $M$
\[
T_{E_1}E_2=\mathcal{H}D_{\mathcal{V}E_1}\mathcal{V}E_2+\mathcal{V}%
D_{\mathcal{V}E_1}\mathcal{H}E_2,\] and
\[A_{E_1}E_2=\mathcal{H}D_{\mathcal{H}E_1}\mathcal{V}E_2+\mathcal{V}%
D_{\mathcal{H}E_1}\mathcal{H}E_2.\]
\end{defn}

In the above definition, $\mathcal{H}$ is the projection on the
horizontal subspace of $TM$, and $\mathcal{V}$ is the projection
on the vertical subspace of $TM$.

The scalar curvature of the Riemannian submersion can be expressed
in terms of the invariants $A$ and $T$ and their covariant
derivatives. Let $(X_i)$ be a local orthonormal basis of
$\mathcal{H}_x$ and $(U_i)$ be a local orthonormal basis of
$\mathcal{V}_x$. The horizontal vector field on $M$ \[ N=\sum_j
T_{U_j}U_j\] is called the mean curvature vector. Denote
\[\check{\grd}N=-\sum_i\sum_j((D_{X_i}T)_{U_j}U_j,X_i),\]
\[|A|^2=\sum_i\sum_j(A_{X_i}U_j,A_{X_i}U_j),\]
\[|T|^2=\sum_i\sum_j(T_{U_j}X_i,T_{U_j}X_i).\]
The following is Corollary 9.37 from~\cite{Bes87}.

\begin{prop} Let $(M,g)$ be a Riemannian submersion over
$(B,\check{g})$ and $(F_b, \hat{g}_b)$ is the fiber over $b\in B$
with the restriction metric $\hat{g}_b=g|_{F_b}$. Let $\grk$,
$\check{\grk}$, $\hat{\grk}$ be the scalar curvatures of the
corresponding metrics. Then \begin{eqnarray}
\grk=\check{\grk}\circ\grp + \hat{\grk} -|A|^2 -|T|^2 -|N|^2
-2\check{\grd}N.\label{warpcurv}\end{eqnarray}
\end{prop}

\subsection{Warped products over one-dimensional base with
spherical fiber}

We now restrict to the case where the base $B$ is one-dimensional
and the fiber is an $(n-1)$-dimensional sphere. Let $B=(0,T)$ and
$S^{n-1}$ is the standard Riemannian sphere with the metric
$d\grx^2$. The object of our investigation is the space of warped
metrics over $B$ with the fiber $S^{n-1}$. Any such metric can be
written as
\[ h=dt^2+f^2d\grx^2.\]

The conditions under which such a metric can be extended to a
smooth metric on an $n$-dimensional disc are given in the
following proposition, which is Lemma 9.114 in~\cite{Bes87}.

\begin{prop}
If we identify $\{x\in\Rb^n|0<|x|<T\}$ with $(0,T)\times S^{n-1}$
in polar coordinates, the smooth Riemannian metric
$dt^2+f(t)^2d\grx^2$ extends to a smooth Riemannian metric on
$\{x\in\Rb^n|\,|x|< T\}$ if and only if $f$ is the restriction to
$(0,T)$ of a smooth odd function on $(-T,T)$ with $f'(0)=1$.
\label{smoothdisc}
\end{prop}

\begin{prop} The scalar curvature $\grk$ for the metric $h$ is
\begin{eqnarray}
\grk=(n-1)\left((n-2)\frac{1-f'^2}{f^2}-2\frac{f''}{f}\right).\label{curveqn}
\end{eqnarray}
\end{prop}

\begin{proof}

It is clear that in this case the tensorial invariant $A$ is $0$.
From Proposition 9.104~\cite{Bes87} it follows that $N$ is basic
and $\grp$-related to the vector field $-(n-1)\frac{f'}{f}$ on
$B$. It follows that
\[ |N|^2=(n-1)^2\frac{f'^2}{f^2}.\]
From the formula 9.105e in~\cite{Bes87} we obtain
\[
\check{\grd}N=(n-1)\left(\frac{f''}{f}+\frac{f'^2}{f^2}\right)\circ\grp.\]
Also, from 9.103 and 9.104 we get $T_{U_j}
X_i=-\frac{(N,X_i)}{(n-1)}U_j$ and
\begin{eqnarray*}
|T|^2&=&\sum_{i,j}(T_{U_j}X_i,T_{U_j}X_i)=\sum_{i,j}\left(\frac{(N,X_i)}{(n-1)}%
U_j,\frac{(N,X_i)}{(n-1)}U_j\right)\\
&=&
\sum_i\frac{1}{(n-1)}(N,X_i)^2=\frac{|N|^2}{(n-1)}=(n-1)\frac{f'^2}{f^2}.
\end{eqnarray*}

Since our base is $1$-dimensional, the scalar curvature
$\check{\grk}=0$. The scalar curvature of $(n-1)$-dimensional
standard Riemannian sphere of radius $f$ is equal to
$\frac{(n-1)(n-2)}{f^2}$. Substituting into~\ref{warpcurv} we
obtain
\begin{eqnarray*}
\grk &=&
\frac{(n-1)(n-2)}{f^2}-2(n-1)\left(\frac{f''}{f}+\frac{f'^2}{f^2}\right)%
-(n-1)^2\frac{f'^2}{f^2}-(n-1)\frac{f'^2}{f^2}\\
&=& (n-1)\left((n-2)\frac{1-f'^2}{f^2}-2\frac{f''}{f}\right).
\end{eqnarray*}
\end{proof}

Any smooth odd function $G$ on $(-T_0,T_0)$ gives rise to a radial
diffeomorphism of euclidean discs
\begin{equation*}
\grUpsilon(G)\colon D_{T_0}\rightarrow D_{G(T_0)},
\end{equation*}
defined by the formula
\begin{equation}\label{diffeofunc}
x\to \frac{G(|x|)}{|x|}x.
\end{equation}

For an even function $g$ on $[-T_0, T_0]$, we set $G(t)=\int_0^t
g(\grt)d\grt$ to obtain a diffeomorphism from $D_{T_0}$ to $D_T$,
where $T=G(T_0)$.

In the proof of Theorem~\ref{locdef}, we will need to deform such
a diffeomorphism to a diffeomorphism which, near the origin, is
multiplication by a constant.

\begin{prop}\label{diffeodeform} Let $\grUpsilon(G)\colon D_{T_0}\rightarrow D_T$ be
the diffeomorphism of euclidean discs, corresponding to an even
function $g$ on $[0,T_0]$, $T:=\int_0^{T_0}g(\grt)d\grt$. Given
any number $\grn<T$, let $\grn^*$ be such that
$\int_0^{\grn^*}g(\grt)d\grt=\grn$. Then there exists a canonical
continuous deformation $\grUpsilon_{\grl}(G)\colon
D_{T_0}\rightarrow D_T$, $\grl\in[0,1]$ such that (i)
$\grUpsilon_0(G)=\grUpsilon(G)$; (ii)
$\grUpsilon_1(G)(x)=\frac{\grn}{\grn^*}x$ on $D_{\grn^*}$; (iii)
$\grUpsilon_{\grl}(G)=\grUpsilon(G)$ on an annulus
$\frac{\grn^*+T_0}{2}\le |x| \le T_0$.
\end{prop}

\begin{proof}
For the proof it is enough to canonically construct a function
$g_1$ such that $g_1(t)=\frac{\grn}{\grn^*}$ on $[0,\grn^*]$,
$g_1=g$ on $[\frac{\grn^*+T_0}{2},T_0]$, and
$\int_0^{T_0}g_1(\grt)d\grt=\int_0^{T_0}g(\grt)d\grt$. Then the
required deformation is given by
\begin{equation*}
\grUpsilon_{\grl}(G)=\grUpsilon((1-\grl)G+\grl G_1).
\end{equation*}
Let
\begin{eqnarray*}
A&=&\int_{\grn^*}^{(\grn^*+T_0)/2}g(\grt)d\grt,\\
A'&=&\min\left(A,\frac{\grn}{2}\left(\frac{T_0}{\grn^*}-1\right)\right).
\end{eqnarray*}
From the Lemma~\ref{cutofflemma} below we can canonically
construct two cutoff functions $\grph_1$ and $\grph_2$ with the
following properties: (i) $\grph_1$ is a decreasing function,
equal to $\frac{\grn}{\grn^*}$ when $0\le t\le \grn^*$, and equal
to $0$ when $\frac{\grn^*+T_0}{2}\le t$; $\grph_2$ is an
increasing function, equal to $0$ on $[0,\grn^*]$, equal to $1$
when $\frac{\grn^*+T_0}{2}\le t$; (iii)
\begin{eqnarray*}
B=\int_{\grn^*}^{(\grn^*+T_0)/2}\grph_1(\grt)d\grt&=& \frac 12 A',\\
\int_{\grn^*}^{(\grn^*+T_0)/2}\grph_2(\grt)g(\grt)d\grt&=&\left(1-
\frac {B}{A}\right) \int_{\grn^*}^{(\grn^*+T_0)/2}g(\grt)d\grt.
\end{eqnarray*}
The function $g_1=\grph_1+\grph_2g$ has all the required
properties. \end{proof}

In general, a warped metric $h$ in the disc $D_{T_0}$ is given by
the formula
\[ h=g^2(t)dt^2+f^2(t)d\grx^2.\] One can think of $g$ as a
diffeomorphism from $D_{T_0}$ to $D_{T}$, $T=G(T_0)$, which
provides an isometry between $h$ and
$\tilde{h}=dt^2+\tilde{f}^2d\grx^2$ with
$\tilde{f}(t)=f(G^{-1}(t))$. The curvature of $h$ is
\begin{eqnarray}
\grk &=&\frac{(n-1)}{f^2g^3}\left((n-2)(g^3-f'^2g)-2f''f
g+2f'fg'\right).\label{curveqn1}
\end{eqnarray}

We fix a family of diffeomorphisms of euclidean discs
\begin{equation}\label{Upsilon}
\grUpsilon_{T_1T_2}\colon D_{T_1}\rightarrow D_{T_2},
\end{equation}
which are radial, i.e. $\grUpsilon_{T_1T_2}(x)=\grl(|x|) x$ with
$\grl > 0$, and are radial isometries near the center and the
boundary of the disc.

This allows us to consider a family of metrics
$h_{\grl}=(g_{\grl},f_{\grl})$, $\grl\in [0,1]$ on discs
$D_{T_{\grl}}$ as metrics on $D_{T_0}$, our fixed size disc.

\subsection{Deformation of warped metrics in the disc}

We deform a warped metric by taking a composition
$f(\grps_{\grl}(t))$, where $\grps_{\grl}$ is a family of
nondecreasing smooth functions. For example, if the family
$\grps_{\grl}$ is such that for some point $c$,
$\grps_1^{(n)}(c)=0$ for all $n\ge 1$, then it means we have
created a collar for the metric $(1,f(\grps_1(t)))$.

Creating a collar with certain properties is the most delicate and
hard step in the deformation process. Lemma~\ref{DefLemma1}
together with the curvature equation~\ref{curveqn} tells us that
for any warped metric we can create an arbitrarily large amount of
curvature arbitrarily close to the center of the disc.
Lemma~\ref{DefLemma2} tells us that once we have this large amount
of curvature, we can bend out a collar. This deformation is
realized in Lemma~\ref{collar}.

In Lemmas~\ref{DefLemma1} and~\ref{DefLemma2} below, we construct
a family of functions defined on varying size intervals. From the
formula~\ref{Upsilon}, it makes sense to talk about continuity of
such deformations with respect to the Fr\'echet topology on
$C^{\infty}[0,T_0]$.

\begin{lem}\label{DefLemma1} Let $C_1\le 1$, $t^*\le T_0/2$, and $\gra\in(0,t^*/2)$ be
positive numbers. Then there exist positive continuous functions
$C(C_1,t^*)$, $T(\grl)$,  $\grl\in[0,1]$, and a family of
functions $\tilde{\grps}_{\grl,\gra}$ on $[0,T(\grl)]$, see
Figure~\ref{Figure31} below, continuously depending on $C_1$,
$\gra$, $t^*$, and $\grl$ with the following properties:
\begin{enumerate}
\item $\tilde{\grps}_{\grl,\gra}(0)=0$,
$\tilde{\grps}_{\grl,\gra}(T(\grl))=t^*$;%
\item $\tilde{\grps}_{0,\gra}(t)=t$;%
\item $\tilde{\grps}_{\grl,\gra}''\le C_1$;%
\item $\tilde{\grps}_{\grl,\gra}''\le 0$ on $\left[\frac{\gra t^*}{40},\frac{\gra t^*}{20}\right]$; %
\item $\tilde{\grps}_{\grl,\gra}''\ge 0$, when  $\tilde{\grps}_{\grl,\gra}(t)%
\in\left[\frac{8}{10}t^*,\frac{9}{10}t^*\right]$;%
\item $\tilde{\grps}_{\grl,\gra}'(t)=1$, for
$\tilde{\grps}_{\grl,\gra}(t)\in
\left[0,\frac{\gra t^*}{40}\right]\cup \left[\frac{9}{10}t^*,t^*\right]$;%
\item $0<1-C(C_1,t^*)\le \tilde{\grps}_{\grl,\gra}'(t)$, all $t\in[0,T(\grl)]$;%
\item on an interval
$\left[\frac{\gra}{10},\frac{8}{10}t^*\right]$,
$\tilde{\grps}_{\grl,\gra}$ follows a straight line\\
with the slope $1-\grl C(C_1,t^*)$.
\end{enumerate}
\end{lem}
\begin{figure}
\centering
\includegraphics{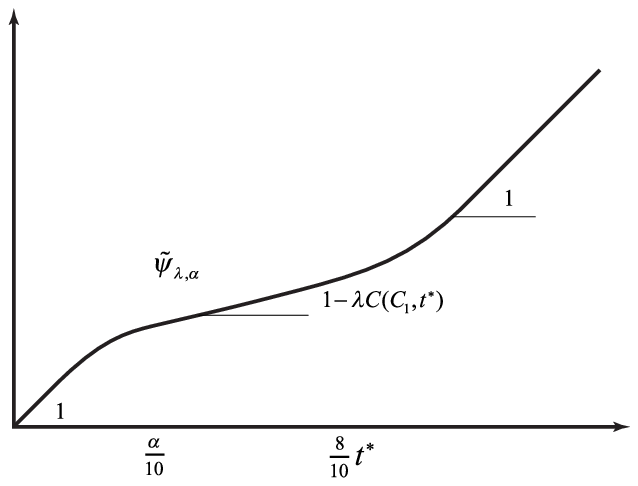}
\caption{}\label{Figure31}
\end{figure}
\begin{proof}

We fix a smooth function $\grf_0$ on $(-\infty,\infty)$, not
identically $0$, such that $\grf_0$ is $0$ on
$(-\infty,0]\cup[1/20,\infty)$, $0\le\grf_0\le 1$. The function
$C$ can now be defined as
\begin{eqnarray}
C(C_1,t^*)=C_1\int_0^{t^*/20}\grf_0\left(\frac{\grs}{t^*}\right)d\grs=C_1
t^*\cdot\mathrm{Const}.
\end{eqnarray}%
The family $\tilde{\grps}_{\grl,\gra}$ will be defined as a
solution to the differential equation
\begin{eqnarray}\tilde{\grps}_{\grl,\gra}''(t)=\grf_{\grl,\gra}(t)\label{diffeqn1}\end{eqnarray}
with the initial conditions $\tilde{\grps}_{\grl,\gra}(0)=0$,
$\tilde{\grps}_{\grl,\gra}'(0)=1$, and a careful choice of the
function $\grf_{\grl,\gra}(t)$ on $[0,T(C_1,t^*)]$. To begin
constructing function $\grf_{\grl,\gra}(t)$, we consider the
equation~\ref{diffeqn1} with the same initial conditions and a
function $\grf_{0,\grl,\gra}(t)=-\frac{2}{\gra}\grl
C_1\grf_0\left(\frac{2t}{\gra t^*}-1/20\right)$ on
$\left[\frac{\gra t^*}{40},\frac{\gra t^*}{20}\right]$ and $0$
otherwise. It is clear that the solution will be a smooth
increasing function whose graph in the $(t,s)$-plane is a curve
that follows a straight line with the slope $1-\grl C(C_1,t^*)>0$
on $t>\frac{\gra t^*}{20}$. It will intercept the horizontal line
$s=\frac{8t^*}{10}$ at some point $t_{\grl}$ which continuously
depends on all possible parameters. We can now define the function
$\grf_{\grl,\gra}$ as follows:
\begin{equation}
\grf_{\grl,\gra}(t)=\begin{cases} -\frac{2}{\gra}\grl
C_1\grf_0\left(\frac{2t}{\gra t^*}-1/20\right) & \textrm{ if } t\in%
\left[\frac{\gra t^*}{40},\frac{\gra t^*}{20}\right]\\
\grl C_1\grf_0\left(\frac{t-t_{\grl}}{t^*}\right)& \textrm{ if }
t\in\left[t_{\grl},t_{\grl}+\frac{t^*}{20}\right]\\
0 & \textrm{ otherwise}\end{cases}.
\end{equation}
Clearly $\grf_{\grl, \gra}$ is a smooth function. One can easily
check that the solution
\begin{equation}
\tilde{\grps}_{\grl,\gra}(t)=t+\int^t_0\int_0^{\grt}\grf_{\grl,\gra}(\grs)d\grs
d\grt
\end{equation}
to the differential equation~\ref{diffeqn1} with this function has
all the required properties. This finishes the proof.\end{proof}

\begin{lem}\label{cutofflemma} Let $\grph_0$ be a smooth nonincreasing cutoff function on
$[0,1]$, i.e. $\grph_0(t)=1$ for $t\le 0$, and $\grph_0(t)=0$ for
$t\ge 1$. Let $0<c< 1$. Then, for any continuous nonnegative $f$
on $[a,b]$, $a\ne b$, which is not identically $0$, there exists a
unique number $\grm
>0$ such that
\begin{equation*}
\int_a^b
\grph_0\left(\left(\frac{t-a}{b-a}\right)^{\grm}\right)f(t)dt=c\int_a^bf(t)dt.
\end{equation*}
The number $\grm$ is a continuous function of $f$ and $c$.
\end{lem}

\begin{proof}
This is a corollary of the Dominated Convergence
theorem.\end{proof}

\begin{figure}

\centering
\includegraphics{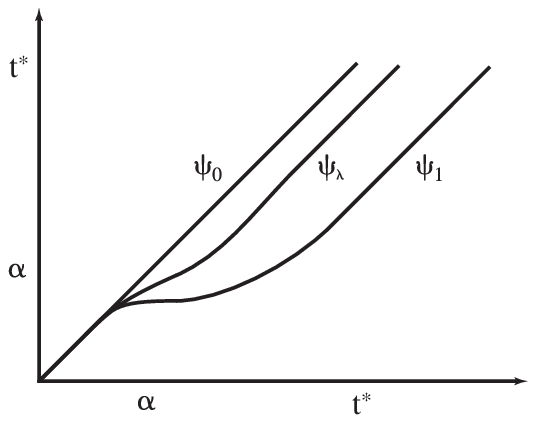}
\caption{}\label{Figure5}

\end{figure}

\begin{lem}\label{DefLemma2} Let $C_1\le 1$ and $t^*\le T_0/2$ be two positive
numbers. Then there exist positive continuous functions
$\gra=\gra(C_1,t^*)$, $T(\grl)$, $\gra\in(0,t^*/2)$ and
$\grl\in[0,1]$ and a family of functions $\grps_{\grl}$ on
$[0,T(\grl)]$, as in Figure~\ref{Figure5} below, continuously
depending on $\grl$ with the following properties:
\begin{enumerate}
\item $\grps_{\grl}(0)=0$,
$\grps_{\grl}(T(\grl))=t^*$;%
\item $\grps_{0}(t)=t$;%
\item $\grps_{\grl}''(t)\le \frac{C_1}{t}$ for all $t>0$;%
\item $\grps_{\grl}'(t)=1$, for
$\grps_{\grl}(t)\in\left[0,\frac{\gra}{10}\right]\cup
\left[\frac{9}{10}t^*,t^*\right]$;%
\item $0\le \grps_{\grl}'\le 1$ ;%
\item $\grps_1^{(n)}\left(\gra\right)=0$, all $n\ge 1$.
\end{enumerate}
\end{lem}

\begin{proof}

We construct $\grps_{\grl}$ as a solution to a second order
ordinary differential equation
\begin{equation}\label{diffeqn2}
\grps_{\grl}''(t)=\grf_{\grl}(t)
\end{equation}
with the initial conditions $\grps_{\grl}(0)=0$ and
$\grps_{\grl}'(0)=1$. The function $\grf_{\grl}$ will be defined
as $\grf_{\grl}=\grl\grf_1$, where $\grf_1$ is constructed in a
few steps. First, we start with the differential equation
\begin{equation*}
f''(t)=\frac{C_1}{t}
\end{equation*}
with the initial conditions $f(\gra)=\gre$, $f'(\gra)=0$ with
$\gre\le\gra$. The solution to this equation is $f(t)=C_1(t\ln t
-t)- (C_1\ln\gra)t+\gra C_1 +\gre$. Let $t_0=\gra
e^{\frac{1.1}{C_1}}$. Then $f'(t_0)=1.1$ and $f(t_0)=\gra\left(1.1
e^{\frac{1.1}{C_1}}+C_1-C_1e^{\frac{1.1}{C_1}}\right)+\gre$. The
function $h(t):=e^{\frac{1}{t}}-te^{\frac{1}{t}}+t$ is a
decreasing positive function on $(0,\infty)$ that goes to $\infty$
as $t$ goes to $0$, and goes to $0$ as $t$ goes to $\infty$. Also,
if $t>0$, then $h(t)< e^{\frac{1}{t}}$. If $\gra$ is small enough
then $f(t_0)$ will be less than any prescribed positive number. We
set $\gra=\frac{9}{20}t^* e^{-\frac{1.1}{C_1}}$. For this choice
of $\gra$ we have that $\gra < \frac{t^*}{2}$, $t_0 <
\frac{9}{10}t^*$, and $f(t_0)<\frac{9}{10}t^*$. Let $\gra_1$ and
$t_1$ be unique points such that
$\int_{\gra}^{\gra_1}\frac{C_1}{\grt}d\grt=\int_{t_1}^{t_0}\frac{C_1}{\grt}d\grt=0.1$.
From the Lemma~\ref{cutofflemma} we can canonically construct two
cutoff functions $\grph_1$ and $\grph_2$ such that for the smooth
function
\begin{equation*}
\tilde{\grf}_1(t)=\begin{cases} \grph_1(t)\frac{C_1}{t}& t\in
[\gra,\gra_1]\\
\frac{C_1}{t}& t\in [\gra_1, t_1]\\
\grph_2(t)\frac{C_1}{t} & t\in [t_1,t_0]\\
0 & \mathrm{ otherwise }
\end{cases}
\end{equation*}
we have $\int_{\gra}^{t_0}\tilde{\grf}_1(\grt)d\grt=1$. We fix a
smooth nonnegative bump function $\tilde{\grf}_2$ which is $0$
outside $[0,1]$ and $\int_{0}^{1}\tilde{\grf}_2(\grt)d\grt=1$ and
define $\grf_1$ as
\begin{equation*}
\grf_1(t)=\begin{cases}
-\frac{10}{9\gra}\tilde{\grf}_2\left(\frac{10t-\gra}{9\gra}\right)
&
t\in\left[\frac{\gra}{10},\gra\right]\\
\tilde{\grf}_1(t) & t\in[\gra, t_0]\\
0 & \mathrm{otherwise} \end{cases}.
\end{equation*}
The solution to the equation~\ref{diffeqn2} is
\begin{equation}
\grps_{\grl}(t)=t+\int_0^t\int_0^{\grt}\grl\grf_1(\grs)d\grs
d\grt.
\end{equation}
On the half interval $[t_0,\infty)$ the function $\grps_{\grl}$ is
strictly increasing for all $\grl$. Therefore, its graph in the
$(t,s)$-plane has a unique intersection $T(\grl)$ with the
horizontal line $s=t^*$.

The family $\grps_{\grl}$ has all the required properties.
Properties $1$, $2$, $3$, and $5$ are obvious from the
construction. The property $6$ follows from the properties of
cutoff functions. The first part of the property $4$ is clear. For
the second part, it is enough to show that $\grps_{\grl}(t_0)\le
\frac{9t^*}{10}$. From the construction we have
\begin{equation*}
\int_{\gra}^{t_0}\int_{\gra}^{\grt}\grf_1(\grs)d\grs d\grt
<\frac{9t^*}{10}-\gra.
\end{equation*}
Then,
\begin{eqnarray*}
\grps_{\grl}(t_0)&\le&
\gra+\int_{\gra}^{t_0}(1-\grl)+\grl\int_{\gra}^{t_0}\int_{\gra}^{\grt}\grf_1(\grs)d\grs
d\grt\\
&<& \gra + (1-\grl)(t_0-\gra)+\grl\left(\frac{9t^*}{10}-\gra\right)\\
&<& \gra +
(1-\grl)\left(\frac{9t^*}{10}-\gra\right)+\grl\left(\frac{9t^*}{10}-\gra\right)=\frac{9t^*}{10}.
\end{eqnarray*}\end{proof}

\subsection{A deformation of $W$ into $W^{\mathrm{loc}}$}

Recall that we defined $W$ to be the space of warped metrics in
the disc $D_{T_0}$ whose scalar curvature is strictly bounded from
below by $B\ge 0$
\begin{eqnarray}
\grk=(n-1)\left((n-2)\frac{1-f'^2}{f^2}-2\frac{f''}{f}\right)>B.\label{curveqnB}
\end{eqnarray}

\begin{defn}\label{loc} Let $h\in W$. We call $h$ a \emph{locally torpedo} metric
if, for some number $c$ with $T_0\ge c >0$, the metric $h$ is a
warped metric $\left(\frac{T_0}{c},
f_0\left(\frac{T_0}{c}t\right)\right)$ in the disc $D_{c}$. Here,
$f_0$ is the warping function of our fixed torpedo metric $g_0$.
The space of all such metrics together with the induced topology
is denoted by $W^{\mathrm{loc}}$.
\end{defn}

\begin{rem}\label{c}
The number $c$ in the above definition is a continuous function on
$W^{\mathrm{loc}}$, $c\colon W^{\mathrm{loc}}\rightarrow (0,T_0]$.
To see this, fix a vector $u_0\in T_{0}D_{T_0}$ which has unit
length with respect to the standard flat metric on $D_{T_0}$.
Then, $c(h)=\frac{T_0}{\sqrt{h(u_0,u_0)}}$. It is easy to see that
$c(h)$ does not depend on the choice of $u_0$.
\end{rem}

The following two lemmas are cornerstones in the construction of a
deformation of $W$ into $W^{\mathrm{loc}}$.

\begin{lem}\label{collar}
There exists a continuous function $\grs\colon W\rightarrow
(0,T_0/2]$ and a continuous map
\[ \grPs_1\colon W\times I\rightarrow W \]
such that $\grPs_1(\cdot,0) = \mathrm{Id}_{W}$ and for the metric
$\tilde{h}=\grPs_1(h,1)=(1,\tilde{f})$ we have for all $n\ge 1$,
the $n$-th derivative $\tilde{f}^{(n)}(\grs(h))=0$.
\end{lem}

\begin{proof}

We write a metric $h\in W$ as a pair $(1,f)$ on $[0,T(h)]$ and
define two continuous functions
\[ \grr_0(t)=\min_{0\le\grt\le t} \frac{f''(\grt)}{2f'''(0)\grt}T(h) \]
and
\[\grr_1(t)=\min_{0\le\grt\le t}\frac{f'(\grt)}{2f'(0)}T(h) \]
be two nonincreasing functions on on $[0,T(h)/2]$, which are equal
to $T(h)/2$ at $0$. It follows that there exist unique numbers
$t_0,\,t_1\in[0,T(h)/2]$ such that $t_0=\grr_0(t_0)$,
$t_1=\grr_1(t_1)$. We define $t^*:=\min(t_0,t_1, T_0/2)$. It
follows that on $[0,t^*]$ the function $f$ satisfies $0<f'(t)<1$
and $f''(t)<0$, when $t>0$ and$f'(0)=1$ and $f''(0)=0$.

Let $B'=\max(\frac 12, \bar{B}')$, where
\[ \bar{B}'=\max_{[0,\frac{9}{10}t^*]}\frac{(n-1)(n-2)%
(1-f'^2)-Bf^2}{2(n-1)ff''},\]%
and
\[ B''=\min_{[\frac
{2\sqrt{2}}{5}t^*,\frac{9}{10}t^*]}\left(\frac{(n-2)}{8}\frac{1-f'^2}{ff'}-%
\frac{f''}{4f'}-\frac{Bf}{8(n-1)f'}\right).\]%
If $C_1=\min (\tilde{C}_1, B'')$, where $\tilde{C}_1$ is such that
$1-C(\tilde{C}_1,t^*)=\sqrt{\frac 12 +\frac{B'}{2}}$, then
$1-C(C_1, t^*)\ge\sqrt{\frac 12 +\frac{B'}{2}}\ge\sqrt{\bar{B}'}$.
For a choice
\[ t^{**}=\min\left(\frac{8}{10}t^*,
\sqrt{(n-1)\frac{1-(1-C(C_1,t^*))^2}{B}}\right),\] %
and
\[C_{1,2}=\min\left(\frac{(n-2)}{8}\left(1-(1-C(C_1,t^*))^2\right),1\right)\]%
we set $\gra=\gra(C_{1,2},t^{**})$.

Let $\tilde{\grps}_{\grl,\gra}$ be a family of curves from
lemma~\ref{DefLemma1} above, corresponding to $C_1$, $t^{*}$, and
$\gra$; $\grps_{\gra}$ a family of curves from
lemma~\ref{DefLemma2} above, corresponding to $C_{1,2}$ and
$t^{**}$. We define
\begin{equation}
\grPs_1(h,\grl)=\begin{cases}(1,f(\tilde{\grps}_{2\grl,\gra}(t))),&
\grl\in[0,1/2]\\
(1,f(\tilde{\grps}_{1,\gra}(\grps_{2(\grl-1/2)}(t)))),&
\grl\in[1/2,1]
\end{cases}.
\end{equation}

We have to show that the scalar curvature of $\grPs_1(h,\grl)$ is
greater than $B$ for all $\grl\in[0,1]$. If the warping function
is $f(\grps(t))$ then the scalar curvature is
\begin{equation}
(n-1)\left((n-2)\frac{1-f'(\grps)^2\grps'^2}{f^2(\grps)}-2\frac{f''(\grps)%
\grps'^2}{f(\grps)}-2\frac{f'(\grps)\grps''}{f(\grps)}\right).
\end{equation}

During the first part of the deformation $\grps''(t)\le 0$ until
$t$ reaches $\frac{8}{10t^*}$. For the scalar curvature to be
greater than $B$ for these values of $t$, it suffices to have
\begin{equation*}
\frac{(n-1)(n-2)(1-f'^2(\grps)\grps'^2)-Bf(\grps)^2}{2(n-1)f(\grps)f''(\grps)}<\grps'^2.
\end{equation*}
This is a consequence of our choice of $C_1$. The only place where
$\grps''$ might be nonnegative is the interval
$[\frac{8}{10}t^*,\frac{9}{10}t^*]$. On this interval it suffices
for $\grps$ to satisfy
\begin{equation}\label{defestimate}
\frac{(n-2)}{2}\frac{(1-f'(\grps)\grps'^2)}{f(\grps)f'(\grps)}-\frac{%
f''(\grps)}{f'(\grps)}\grps'^2-\frac{B}{2(n-1)}\frac{f(\grps)}{f'(\grps)}%
>\grps''
\end{equation}
for the scalar curvature to be greater than $B$. However, on this
interval $\frac{9}{10}t^*\ge \grps(t)\ge \frac{2\sqrt{2}}{5}t^*$
and $\grps'^2\ge \left(\frac 12 + \frac{B'}{2}\right)$. Therefore,
\begin{eqnarray*}
\lefteqn{\frac{(n-2)}{2}\frac{(1-f'^2(\grps))}{f(\grps)f'(\grps)}-
\frac{f''(\grps)}{f'(\grps)}\left(\frac 12 +\frac{B'}{2}\right)
-\frac{B}{2(n-1)}\frac{f(\grps)}{f'(\grps)}}\\
&>& \frac{(n-2)}{2}\frac{(1-f'^2(\grps))}{f(\grps)f'(\grps)}%
-\frac{f''(\grps)}{2f'(\grps)}-\frac{(n-1)(n-2)(1-f'^2(\grps))-B
f^2(\grps)}{4(n-1)f(\grps)f'(\grps)}\\
&-&\frac{B}{2(n-1)}\frac{f(\grps)}{f'(\grps)}\\
&>&\frac{(n-2)(1-f'^2(\grps))}{4f(\grps)f'(\grps)}-\frac{f''(\grps)}{2f'(\grps)}%
-\frac{B}{4(n-1)}\frac{f(\grps)}{f'(\grps)}\\
&>&B''\ge\grps''.
\end{eqnarray*}
This shows that the scalar curvature remains strictly bounded from
below by $B$ during the first part of the deformation.

To estimate the scalar curvature during the second deformation it
is enough to look at the inequality~\ref{defestimate}. This is,
because after the first deformation on the interval
$\left[\frac{\gra}{10},t^{**}\right]$ the inequalities $0<f'\le
\tilde{C}=(1-C(C_1,t^*))<1$ will be satisfied. We have
\begin{eqnarray*}
\lefteqn{\frac{(n-2)}{2}\frac{(1-f'^2(\grps)\grps'^2)}{f(\grps)f'(\grps)}-%
\frac{B}{2(n-1)}\frac{f(\grps)}{f'(\grps)}}\\
&\ge&
\frac{(n-2)}{2}\frac{(1-\tilde{C}^2)-\frac{B}{2(n-1)}f^2(\grps)}%
{f(\grps)}\\
&\ge&\frac{(n-2)}{2}\frac{(1-\tilde{C}^2)-\frac{B}{2(n-1)}t^2}{t}\\
&\ge& 2\frac{C_{1,2}}{t}>\grps''(t).
\end{eqnarray*}
We define $\grs(h):=\gra$. Let $(1,\tilde{f}):=\grPs(h,1)$. Then
the property $\tilde{f}^{(n)}(\grs(h))=0$ follows from the
construction and the properties of cutoff functions.\end{proof}

\begin{defn} Let $D_{[T_0,T_1]}^n:=\{x\in\Rb^k\ |\ 0<T_0\le |x|\le T_1\}$
be a euclidean annulus and $h$ is a metric on it. We call $h$ an
\emph{annular} metric if (i) $h$ is warped with the fibre metric
$d\grx^2$ the standard euclidean sphere of radius $1$; (ii) near
the boundary components $S_{T_i}^{n-1}$, $i=0, 1$, we have that
$h=dt^2+r_i^2 d\grx^2$ for some constants $r_i>0$.
\end{defn}

We fix a family of annular metrics $\hat{h}$ which is continuous
in $T_i$, $r_i$ and the scalar curvature of $\hat{h}$ is strictly
bounded from below by $B\ge 0$. This definition makes sense since
one can think of an annular metric as a warped metric on
$\Rb\times S^{n-1}$.

This family $\hat{h}$ allows us to consider a deformation defined
on a subdisc (subannulus) as a deformation in the whole disc
(annulus) that is constant on some neighborhood of the boundary.

\begin{lem}\label{locd} Let $h=(1,f)$ be a warped metric on the disc
$D_{T(h)}$ and $\grs\in(0,T(h))$, $\grs<T_0/2$, be such that
$f^{(n)}(\grs)=0$ for all $n\ge 1$; $0\le f'\le 1$, $f''\le 0$ on
$[0,\grs]$. Then there exists a positive number $\grd_0$ and a
family of metrics $\grPs_2(h,{\grl},\grs)=(g_{\grl},f_{\grl})$,
$\grl\in[0,1]$ on the disc $D_{T_{\grl}}$ such that
$\grPs_2(h,0,\grs)=h$, $g_1=1$, $f_1=f_0$ on $[0,T_0]$, and the
restriction of $h$ to an annulus $D_{[T(h)-\grd_0,T(h)]}$ is
isometric to the restriction of $\grPs_2(h,{\grl},\grs)$ to an
annulus $D_{[T_{\grl}-\grd_0,T_{\grl}]}$ for all $\grl$. The
family $\grPs_2(h,{\grl},\grs)$ is continuous in $h$, $\grs$, and
$\grl$.

\end{lem}

\begin{proof}
We take a linear deformation $f_{\grl}=(1-\grl)f+\grl f_0$ on the
interval $[0,\grs]$. In general, the condition in the
formula~\ref{curveqnB} above is not convex. However, it is convex
when $B=0$ and $f$ is such that $0\ge f'(t)\le 1$, $f''(t)<0$ for
all $f\in[0,\grs]$ ($f_0$ can always be chosen to satisfy these
conditions). That is, the scalar curvature for $(1,f_{\grl})$ is
positive for all $\grl\in [0,1]$, if the scalar curvature for $f$
were positive. Since the family $h_{\grl}=(1,f_{\grl})$ is
compact, it is possible to find a number $1\ge\grn>0$ continuously
depending on $f$ such that the scalar curvature of the family
$\grn h_{\grl}$ in the disc $D_{\grs}$ is bounded from below by
$B$. We then plug this family in the disc $D_{T_0}$ using the
fixed family of annular metrics $\hat{h}$. The existence of
$\grd_0$ follows from the compactness of $[0,1]$.\end{proof}

\begin{thm}\label{locdef} There exist a number $\grd>0$ and a
continuous map
\begin{equation}
\grPs\colon W\times I\rightarrow W,
\end{equation}
such that: (i) $\grPs(h,0)=h$; (ii) $\grPs(h,1)\in
W^{\mathrm{loc}}$; (iii) $\grPs(h,t)=h $ on the annulus
$D_{[T_0-\grd,T_0]}$ for all $t\in[0,1]$.
\end{thm}

\begin{proof}

Let $W'=\grPs_1(W\times 1)$ be the space from the deformation of
lemma~\ref{collar}. A metric $h=(\bar{g},\bar{f})\in W'$ has the
following property. After reparametrization, we can write
$h=(1,f)$. Then the warping function $f$ on $[0, T(h)]$ is an odd
smooth function, $0\le f'\le 1$, $f''\le 0$ on the interval
$[0,\grm(h)]$, where $\grm(h)$ is a continuous function of metric
($\grm$ is the composition of $\grs$ and the diffeomorphism
function of $\bar{g}$) and $f^{(n)}(\grm(h))=0$, all $n\ge 1$. We
also notice from the proof of Lemma~\ref{collar}, that $\grPs_1$
does not change the metric on the collar of the boundary of
$D^n_{T_0}$ of some fixed size $\grd_1$.

The required deformation is defined as
\begin{equation*}
\grPs(h,t):=\begin{cases}\grPs_1(h,2t),& t\in[0,1/2]\\
\grPs_2(\grPs_1(h,1),2t-1,\grm(\grPs_1(h,1))),& t\in[1/2,1]
\end{cases},
\end{equation*}
where $\grPs_2$ is the deformation from Lemma~\ref{locd}. We take
$\grd:=\min({\grd_1,\grd_0})$, where $\grd_0$ is from
Lemma~\ref{locd}.

After Proposition~\ref{diffeodeform}, we may assume that
$\grPs(h,1)\in W^{\mathrm{loc}}$. This completes the proof.
\end{proof}

We finish this section with a lemma that we will need in the proof
of Theorem~\ref{weakdef}.

\begin{lem}\label{annulus}
Let $D_{[t_1,t_2]}$ be a euclidean annulus, and $h$ be a warped
metric on this annulus, i.e., in polar coordinates, $h=g(t)^2
dt^2+ f(t)^2d\grx^2$, where $d\grx^2$ is the standard metric on
$S^{k-1}$. Let $B\ge 0$ be a lower bound for the scalar curvature
of $h$ in $D_{[t_1,t_2]}$. Suppose that for some $0<\gre
<(t_2-t_1)/2$ and some $0<r_0$ we have $g\equiv 1$ and $f\equiv
r_0$ on $[t_1,t_1+\gre]\cup [t_2-\gre,t_2]$. Then there exists a
continuous family of warped metrics $h_{\grl}$, $\grl\in[0,1]$
such that:
\begin{enumerate}
\item[(i)] $h_0=h$;%
\item[(ii)] $h_1=dt^2+r_0^2d\grx^2$;%
\item[(iii)] $\exists\, \grd>0$ such that $h_{\grl}\equiv h$
whenever $|x|\in[t_1,t_1+\grd]\cup[t_2-\grd,t_2]$;%
\item[(iv)] $\grk_{h_{\grl}}(x)> B$, all $x\in D_{[t_1,t_2]}$,
$\grl\in[0,1]$.
\end{enumerate}
\end{lem}

\begin{proof}

(i) The case $B=0$. We may assume $g\equiv 1$.  Set
$f_{\grl}:=(1-\grl)f+\grl r_0$, $\grl\in[0,1]$. From
equation~\ref{curveqn1} there exists $A\ge 1$ such that the scalar
curvature of the family $A^2 dt^2+f_{\grl}^2(t)d\grx^2$ is
positive. One can easily write down the deformation $h_{\grl}$.

(ii) The case $B>0$. We introduce a function \[S(h):=\min_{x\in
D}\grk_h(x),\] where $D$ is the annulus on which $h$ is defined.
For the family $h_{\grl}$ defined above the condition
$S(h_{\grl})>B$ does not have to hold. However, the condition
$S(h_{\grl})>0$ (which was satisfied initially, i.e. when
$\grl=0$) will still hold. This allows to rescale the whole family
$h_{\grl}$ with a single scaling factor $\grn$, so that the scalar
curvature of each metric in the rescaled family is greater than
$B$. We plug this family $\grn h_{\grl}$ into the annulus using
the fixed family of annular metrics $\hat{h}$.\end{proof}

\section{Proof of the main theorem}

Let
\begin{equation*}
B=\max\left(-\min_{x\in N}\grk_{g_N}(x),0\right)
\end{equation*}
and
\begin{equation*}
W(N, \grt)=\left\{g\in\pscm{M^n}|\ \grt^*(g)=g_N+g_w\right\},
\end{equation*}
where $g_w$ is a warped metric in the disc $D^k_{T_0}$ such that
the scalar curvature of $g_w$ is greater than $B$. Analogously, we
define the subspace $W^{\mathrm{loc}}(N,\grt)\subset W(N,\grt)$,
cf. Definition~\ref{loc}.

The following theorem holds.
\begin{thm}\label{weakdef} The subspace $\pscmt{M}$ is a weak deformation retract
of $W(N,\grt)$.
\end{thm}
\begin{proof}
Without loss of generality we may assume that the tubular map
$\grt$ is an embedding on $N\times D^k_{T_1}$, where $T_1>T_0$.

We fix a continuous family $\grUpsilon_{\grl,T}$, with
$(\grl,T)\in[0,1]\times(0,T_0]$, of radial diffeomorphisms of
$D_{T_1}$ with the following properties:
\begin{enumerate}
\item[(i)] $\grUpsilon_{0,T}=\mathrm{Id}$, all $T\in(0,T_0]$;%
\item[(ii)] $\grUpsilon_{1,T}$ acts by multiplication by
$\frac{T}{T_0}$ on $D_{T_0}$;%%
\item[(iii)] if $T^*\ge T_0$ is such that
$\grUpsilon_{1,T}(S_{T^*})=S_{T_0}$, then $\grUpsilon_{1,T}$ is a
radial isometry in some small radial neighborhood of $S_{T^*}$.
\end{enumerate}

From Remark~\ref{c} we have a continuous map \[c\colon
W^{\mathrm{loc}}(N,\grt)\rightarrow (0,T_0].\] The retraction map
$r\colon W(N,\grt)\rightarrow \pscmt{M}$ can now be defined as a
composition
\begin{equation}
W(N,\grt)\xrightarrow{\grPs_1}W^{\mathrm{loc}}(N,\grt)\xrightarrow%
{\grUpsilon}\pscmt{M},
\end{equation}
where $\grPs_1=\grPs(\cdot,1)$ is  from the deformation
$\grPs\colon W\times I\rightarrow W$ that was constructed in
Theorem~\ref{locdef}, and
$\grUpsilon(h)=\grUpsilon_{1,c(h)}^{*}(h)$.

The deformation $D\colon W(N,\grt)\times I\rightarrow W(N,\grt)$
from $\mathrm{Id}$ to $r$ is given as
\begin{equation*}
D(h,\grl)=\begin{cases}\grPs(h,2\grl)&\grl\in[0,1/2]\\
\grUpsilon_{(2\grl-1),c(h')}^*(h')& \grl\in[1/2,1]
\end{cases},
\end{equation*}
where $h'=\grPs(h,1)$. This proves that $i\circ
r\simeq\mathrm{Id}$, i.e. that the map $r$ is the right homotopy
inverse for the natural inclusion map $i$.

To show that $r\circ i\simeq\mathrm{Id}$ on $\pscmt{M}$ we begin
with the following observation. For all $h\in\pscmt{M}$ the
function $c(\grPs(h,1))=T=\mathrm{const}\le T_0/2$. This follows
from the construction of $\grPs$ (see the proofs of
Lemma~\ref{collar} and Theorem~\ref{locdef} for details). This
implies that for all such metrics there exists a constant
$T^*>T_0$ such that $r(h)$ is equal to the torpedo metric $g_0$ on
$D_{T_0}$, is equal to a fixed warped metric
$h_w=g(t)^2dt^2+f(t)^2d\grx^2$ on the annulus $D_{[T_0,T^*]}$ with
$g(t)\equiv 1$ in some small left neighborhood of $T^*$, and is
equal to a pullback metric on the annulus $D_{[T^*,T_1]}$. The
deformation that we are seeking consists of (a) applying
Proposition~\ref{annulus} to deform $h_w$ to $dt^2+r_0^2 d\grx^2$;
(b) contracting the annulus $T_0\le |x|\le T^*$ to the sphere
$S_{T_0}$ and at the same time stretching the annulus $D_{[T^*,
T_1]}$ back to the identity on $D_{[T_0,T_1]}$ by some family of
diffeomorphisms which are radial isometries near the boundary.
\end{proof}

Let $g_s\in\pscm{M}$, $s\in S$, be a continuous family of metrics
and $S$ is a compact manifold. Then from Theorems~\ref{isolemma},
~\ref{alphatwo}, and~\ref{alphathree} we conclude that there
exists a continuous map \[\GL\colon S\times I\rightarrow\pscm{M}\]
such that (i) $\GL(s,0)=g_s$; (ii) $\GL(s,1)\in\pscmt{M}$; (iii)
$\GL(s,0)\in W(N,\grt)$ implies that $\GL(s,t)\in W(N,\grt)$ for
all $t\in[0,1]$. The last property follows from the construction
of the map $\GL$, the definition of the space $W(N,\grt)$, and the
fact that we can always choose $\GL$ to take place inside the
image of the tubular map $\grt$. This deformation is by no means
unique and we call it a $\GL$ deformation for the family $g_s$.

\begin{proof}[Proof of Theorem~\ref{main}]
It is enough to show that all the relative homotopy sets
\begin{equation}\label{hgroups}
\grp_r\left(\pscm{M^n},W(N,\grt)\right)=0, \quad r\ge 1,
\end{equation}
and the inclusion map $i\colon W(N,\grt)\rightarrow \pscm{M}$ is a
bijection between path connected components of both spaces.

We first show that $i$ induces a bijection between the path
components. Let $h_1,\ h_2\in W(N,\grt)$ and $\gra\colon
I\rightarrow\pscm{M^n}$ be a path between $h_1$ and $h_2$. Take
some $\GL\colon I\times I\rightarrow\pscm{M}$ for the family
$\gra$. Then $\GL(t,1)\in\pscmt{M}\subset W(N,\grt)$ and
$\GL(h_1,t),\ \GL(h_2,t)\in W(N,\grt)$ for all $t\in[0,1]$. The
path
\begin{equation*}
\gra'(t):=\begin{cases}\GL(h_1,3t)& t\in[0,1/3]\\
\GL(\gra(3(t-1/3)),1)& t\in[1/3,2/3]\\
\GL(h_2,3(1-t))&t\in[2/3,1]
\end{cases}
\end{equation*}
is the required path.

Let $[\gra]\in\grp_r\left(\pscm{M^n},W(N,\grt)\right)$, where
$\gra\colon (D^r,S^{r-1})\rightarrow(\pscm{M^n},W(N,\grt))$ is a
continuous map. Take a deformation
\[
\GL\colon D^r\times I\rightarrow\pscm{M^n}
\]
for the family $\gra$ to conclude that $[\gra]=0$.\end{proof}

\begin{proof}[Proof of Theorem~\ref{surgery}]
We fix the standard euclidean metrics of radius $1$ on $S^k$ and
on $S^{n-k-1}$ and torpedo metrics on $D^{n-k}$ and on $D^{k+1}$
of size $T_0$ and radius $1$. From the definition of a torpedo
metric there exist a number $\gre>0$ such that, near the boundary
of the disc, both torpedo metrics are products of the standard
sphere with an interval of length $\gre$.

Let $M_0:= M_1-(S^k\times\stackrel{\circ}{D}
{}^{n-k}_{T_0-\gre/2})$. Then
\[M_2=M_0\bigcup_{S^k\times S^{n-k-1}}D^{k+1}_{T_0-\gre/2}\times S^{n-k-1}.\]
This gives us an embedding of $S^{n-k-1}\rightarrow M^n_2$ and a
tubular neighborhood \[\grt_2\colon S^{n-k-1}\times
D^{k+1}_{T_0}\rightarrow M_2^n.\]

From Theorem~\ref{main}, we have $\pscm{M_1}\simeq\pscmt{M_1}$ and
$\pscm{M_2}\simeq\pscmt{M_2}$. But the space $\pscmt{M_1}$ is
homeomorphic to $\pscmt{M_2}$ as can be seen by removing
$S^k\times D^{n-k}_{T_0-\gre/2}$ from $M_1$ and removing
$S^{n-k-1}\times D^{k+1}_{T_0-\gre/2}$ from $M_2$. In both cases
we obtain the manifold $M_0$ with the fixed collar $S^k\times
S^{n-k-1}\times[0,\gre/2]\rightarrow M_0$ and the space of
positive scalar curvature metrics on $M_0$ that restrict to the
fixed product metric on this collar.\end{proof}

\bibliographystyle{amsalpha}
\providecommand{\bysame}{\leavevmode\hbox
to3em{\hrulefill}\thinspace}
\providecommand{\MR}{\relax\ifhmode\unskip\space\fi MR }
\providecommand{\MRhref}[2]{%
  \href{http://www.ams.org/mathscinet-getitem?mr=#1}{#2}
} \providecommand{\href}[2]{#2}

\end{document}